\documentstyle[12pt,twoside, epsf]{article}

\def\picill#1by#2(#3)
{\vbox to #2
{\hrule width #1 height 0pt depth 0pt
\vfill\epsffile{#3}}}

\let \ttorg \tt \def \tt{\ttorg \obeyspaces}

\begin{document}

\date{}

\title{\bf Fibonacci Rectangles}

\author{Louis H. Kauffman \\
Department of Mathematics, Statistics and
Computer Science\\
University of Illinois at Chicago\\
851 South Morgan St., Chicago IL 60607-7045, USA\\
kauffman@uic.edu
}

 \maketitle
  
 \thispagestyle{empty}

\section{Introduction}
 
In Figure 1 we have a rectangle of length $377$ units and width $233$ units.
It is paved with squares of sizes $233 \times 233, 144 \times 144, 89 \times 89, 55 \times 55, 34 \times 34,
21 \times 21, 13 \times 13, 8 \times 8, 5 \times 5, 3 \times 3, 2 \times 2, 1 \times 1, 1 \times 1.$
This attests to the fact that for the Fibonacci numbers

$$1,1,2,3,5,8,13,21,34,55,89,144,233,377, \cdots$$
the sum of the squares of the Fibonacci numbrers up to a given Fibonacci number is equal to the product of that Fibonacci 
number with its successor. For example,
$$ 1^2 + 1^2 + 2^2 + 3^2 + 5^2 + 8^2 = 8 \times 13.$$
\bigbreak

Letting
$$f(0) = 1$$
$$f(1) = 1$$
$$f(2) = 2$$
$$f(3) = 3$$
$$f(4) = 5$$
and 
$$f(n+1) = f(n) + f(n-1),$$
we have that a rectangle of size $f(n) \times f(n+1)$ can be paved with squares of different sizes, all except for a repetition
of the size $1 \times 1$ at the very end of the spiralling process of cutting off squares that are the width of the given
rectangle. 
\bigbreak

It is well-known that this process of cutting off squares can be continued to infinity if we start with
a rectangle that is of the size $\Phi \times 1$ where $\Phi$ is the Golden Ratio $$\Phi = (1 + \sqrt{5})/2).$$ This is no mystery.
Such a process will work when the new rectangle is similar to the original one. That condition is embodied in the equation
$$W/(L-W) = L/W,$$ and taking $W=1,$ we find $1/(L-1) = L,$ whence $$L^2 - L -1 = 0,$$ whose positive root is the Golden Ratio.
\bigbreak

It is also 
well-known that $\Phi$ is the limit of successive ratios of Fibonacci numbers with 
$$ 1 < 3/2 < 8/5 < 21/13 < \cdots < \Phi < \cdots < 13/8 < 5/3 < 2.$$

\begin{center}
$$ \picill5inby3.5in(FibonacciRectangle)  $$
{ \bf Figure 1 -- The Fibonacci Rectangles} \end{center}

We ask: {\bf Is there any other proportion for a rectangle, other than the Golden Proportion, that will allow the process of 
cutting off successive squares to produce an infinite paving of the original rectangle by squares of different sizes?}
\bigbreak

The answer is: {\bf No. The only proportion that allows this pattern is the Golden Ratio!}

\begin{center}
$$ \picill5inby3.5in(Fibonacci2)  $$
{ \bf Figure 2 -- Characterizing the Golden Ratio} \end{center}  

View Figure 2.
Suppose the original rectangle has width $W$ and length $L.$
In order for the process of cutting off the square to produce a new rectangle whose own square cut-off is smaller than
the first square, we need the new rectangle to have width $L-W$ and length $W$ with $ L-W < W.$ 
With this inequality, the square cut off from the new rectangle will be smaller than the first square cut off from the original 
rectangle.
\bigbreak

In order for this pattern of cutting to persist forever, we need an infinite sequence of inequalities, each derived from the
previous one:
$$W < L$$
$$L - W < W$$
$$W - (L - W) < L - W$$
$$L- W - (W - (L - W)) < W - (L - W)$$
ad infinitum.\bigbreak

You can easily see the pattern:
$$W < L$$
$$L - W < W$$
$$2W - L < L - W$$
$$2L - 3W < 2W - L$$
$$5W - 3L < 2L - 3W$$
$$5L - 8W < 5W - 3L$$
$$13W -8L < 5L - 8W$$
ad infinitum.
But now the Fibonacci pattern is apparent! Look at what each of these inequalities says:
$$W < L$$
$$L < 2W$$
$$3W < 2L$$
$$3L < 5W$$
$$8W < 5L$$
$$8L < 13W$$
$$21W < 13L$$
ad infinitum.
These in turn say:
$$L/W > 1$$
$$L/W < 2$$
$$L/W > 3/2$$
$$L/W < 5/3$$
$$L/W > 8/5$$
$$L/W < 13/8$$
$$L/W > 21/13$$
ad infinitum.
\bigbreak

Thus we see that the pattern is that $L/W$ is sandwiched between ratios of Fibonacci numbers just as we know the Golden Ratio is 
sandwiched:
$$ 1 < 3/2 < 8/5 < 21/13 < \cdots < L/W < \cdots < 13/8 < 5/3 < 2,$$
and this implies that $$L/W = \Phi = (1 + \sqrt{5})/2.$$

\section{Matrices}
Here is a matrix formulation of the results in the previous section.
Let 
$$M = \left( \begin{array}{cc}
-1 & 1 \\
 1 & 0 \\
\end{array} \right).$$ Then
$$M \left( \begin{array}{c}
W\\
L\\ \end{array} \right) = 
\left( \begin{array}{c}
L- W\\
W \\ \end{array} \right) =
\left( \begin{array}{c}
W'\\
L'\\ \end{array} \right)$$
where $W'$ and $L'$ are the width and length for the rectangle obtained by cutting off a square of side $W$ from the 
rectangle of size $W \times L.$

Letting
$$\left( \begin{array}{c}
W^{(n)}\\
L^{(n)}\\ \end{array} \right) =
M^{n}\left( \begin{array}{c}
W\\
L\\ \end{array} \right),$$ our demand on the proportions of these rectangles is simply that $0 < W^{(n)} < L^{(n)}$ for all
$n = 1,2,3, \cdots.$

It is easy to see by induction that 
$$M^{n} = (-1)^{n}\left( \begin{array}{cc}
f(n) & -f(n-1) \\
 -f(n-1) & f(n-2)\\
\end{array} \right),$$ where $f(n)$ denotes the Fibonacci series as defined in the first section.
Thus
$$(-1)^{n}\left( \begin{array}{c}
W^{(n)}\\
L^{(n)}\\ \end{array} \right) = \left( \begin{array}{cc}
f(n) & -f(n-1) \\
 -f(n-1) & f(n-2)\\
\end{array} \right)\left( \begin{array}{c}
W\\
L\\ \end{array} \right) =
\left( \begin{array}{c}
f(n)W - f(n-1)L\\
-f(n-1)W + f(n-2)L\\ \end{array} \right).$$
Hence $W^{(n)} < L^{(n)}$ if and only if 
$$f(n)W - f(n-1)L < -f(n-1)W + f(n-2)L$$ for $n$ even, and 
$$f(n)W - f(n-1)L > -f(n-1)W + f(n-2)L$$ for $n$ odd. 
This implies that
$$f(n+1)W < f(n)L$$ for $n$ even, and 
$$f(n+1)W > f(n)L$$ for $n$ odd. Whence
$$L/W > f(n+1)/f(n)$$ for $n$ even, and 
$$L/W < f(n+1)/f(n)$$ for $n$ odd. 
Thus this matrix formulation shows systematically that 
$L/W$ is equal to the Golden Ratio.
\bigbreak

\section {Finitely Squared Rectangles and Other Diversions}
It is a famous problem to pave a rectangle with a finite number of squares of different sizes.
For a history of this problem and its solution using graph theory and the properties of electrical circuits,
the reader should consult the work of W. T. Tutte and his co-authors \cite{Tutte1,Tutte2}.
\bigbreak

Here we have considered one specific problem in the paving of rectangles by an infinite number of distinct squares.
Clearly there is more to be done in this field.
\bigbreak

There are relationships with the theory of knots lurking near this paper. For one thing, if you take a strip of paper and tie it
into a simple overhand knot, and pull it up tight and flat, you will produce a perfect regular pentagon \cite{BANDS}. 
See Figure 3.  The pentagon
is a close relative of the golden ratio, embodying that ratio in the proportion of its chords to its sides. Then again, there are
electrical  circuits in Tutte's work, and there are electrical circuits in back of the knots \cite{KNOTS,GOLD}.

\begin{center}
$$ \picill5inby1.5in(PentTrefoil)  $$
{ \bf Figure 3 -- Trefoil Implicates Pentagon} \end{center}

\noindent {\bf Acknowledgement.} It gives the author great pleasure to thank Bob and Ellen Kaplan for their lecture
(part of the Perimeter Institute Public Lecture Series) on April 7, 2004. The characterization of the Golden Ratio for 
infinite  (successive square cut-off) pavings of rectangles with unequal squares was inspired by the open atmosphere of
questioning and creating, fostered by the Kaplan lecture, and proved by the author during that lecture. This is better than
falling asleep! The author would also like to thank Blanches Descartes for many helpful conversations, as it has been his
privilege to share her office this year (2003 - 2004) while visiting the University of Waterloo and the Perimeter Institute in
Waterloo, Canada.
\bigbreak

\end{document}